\documentclass[oneside,12pt]{amsart}
\usepackage{amssymb, amscd}

\theoremstyle{plain}
\newtheorem{theorem}{Theorem}[section]
\newtheorem{lemma}[theorem]{Lemma}

\newtheorem*{theorem*}{}
\newtheorem{proposition}[theorem]{Proposition}
\newtheorem{corollary}[theorem]{Corollary}

\theoremstyle{definition}

\newtheorem{definition}[theorem]{Definition}

\newtheorem{example}[theorem]{Example}
\newtheorem{remark}[theorem]{Remark}

\def\U{{\bf U}}

\begin{document}

\title{Twisted Orbifold K-theory}

\author{Alejandro Adem and Yongbin Ruan}

\date{\today}
\address{Mathematics Department\\
         University of Wisconsin\\
         Madison, Wisconsin, 53706}
\email{adem@math.wisc.edu, ruan@math.wisc.edu}
\thanks{Both authors were partially supported by the NSF}

\begin{abstract}
We use equivariant
methods to define and study the \emph{orbifold K-theory} of
an orbifold $X$.
Adapting 
techniques from equivariant K-theory, 
we construct a Chern character
and exhibit a multiplicative decomposition for
$K^*_{orb}(X)\otimes\mathbb Q$, in particular showing that it is
additively isomorphic to the orbifold cohomology of $X$. 
A number of examples are
provided. We then use the theory of projective representations to
define the notion of twisted orbifold K--theory in the presence of
discrete torsion. An explicit
expression for this is obtained 
in the case of a global quotient.
\end{abstract}

\maketitle 
\tableofcontents
\section{Motivation}

    During the last twenty years, one major trend has been
    the constant flow of physical ideas into mathematics. Although
    many of them have had a significant
    impact, there are still
    many others which deserve more attention from mathematicians.
    Orbifold string theory is such an example:
    it has been around in physics since 1985; compared to its
    popularity in physics, its influence in mathematics is
    minimal. This project represents an effort to change this picture.

    First, let us give a brief introduction to orbifold string
    theory aiming to explain
    its mathematical content. In the end, we hope
    that it will become clear why we have to study
    orbifold
    K-theory. In 1985, Dixon-Harvey-Vafa-Witten \cite{DHVW} discovered that on
    a singular manifold
    such as an orbifold one can
    also build a smooth consistent string theory.  One of the remarkable
    insights of orbifold string theory
    is that its consistency demands the
    introduction of so called twisted sectors. In another words,
    ordinary cohomology is the \textbf{wrong}
    theory for orbifolds; the
    correct one seems to be {\em orbifold cohomology},
    which has been
    constructed by Chen-Ruan \cite{CR}. Orbifold cohomology
    (not ordinary cohomology) should
    fulfill the role of cohomology for
    smooth manifolds. Moreover, on smooth manifolds, orbifold
    cohomology is the same as ordinary cohomology.

    One of the original motivations of orbifold string theory is to
    capture the information of its {\sl crepant resolution}
    (see \cite{Jo} page 126 for definitions). In fact, it
    is generally believed that orbifold string theory is
    equivalent to the ordinary string theory of its crepant resolution
    in some fashion.  More precisely, one can associate a super
    conformal field theory to each 3-dimensional Calabi-Yau orbifold or smooth
    manifold. The orbifold super conformal field theory and ordinary
    super conformal field theory of its crepant resolution appear
    to be two members of a family of superconformal field theories.
    To explore its mathematical implications, it is useful to
    consider its low energy effective theory. Low energy effective
    theory is parameterized by the moduli space of vacua. At each
    vacuum, one can associate the orbifold cohomology (for the orbifold
    case) or ordinary cohomology (for the smooth case). Then there is
    a natural metric on orbifold cohomology corresponding to
    Poincar\'e duality. After quantization, the metric is corrected by
     the quantum corrections  given by
    orbifold quantum cohomology (see \cite{CR1}).

    To compare
    the theory at an orbifold point with that of its crepant
    resolution, one has to consider the variation of low energy
    effective theory  as we vary the point of moduli space of
    superconformal field theory. It is known that the orbifold points are often
    singularities. Its  monodromy could
    interchange $H^{1,1}$ and $H^{2,2}$. Therefore, there is no
    hope to preserve the grading of orbifold cohomology. However,
    it still preserves the parity. This naturally leads to
    the concept of \emph{orbifold K-theory}
    \footnote{Throughout this paper we shall be dealing
    exclusively with (equivariant) complex $K$--theory, i.e. the
    generalized cohomology
    theory arising from complex vector bundles} which we shall define
    in section 4.
    In this context, a natural
    conjecture arising from orbifold string theory is
   \medskip

    \noindent
    {\bf K-Orbifold String Theory Conjecture:}
    {\it There is a natural additive isomorphism between
    orbifold K-theory and the ordinary K-theory of its crepant
    resolution.}
    \vskip 0.1in

    Another physical motivation comes from the recent discovery in physics that
    the
    D-brane charge of an ordinary string theory is described by the K-theory
    of its underlying smooth manifold.
    For orbifold string
    theory, it is obvious that the orbifold D-brane charge should be
    described by our orbifold K-theory.

    An obvious definition of orbifold K-theory is the
    K-theory associated to
    orbifold vector bundles. However, to be a correct definition
    for our purposes, we have to obtain at least an additive
    isomorphism to orbifold
    cohomology. It is well-known that any reduced orbifold can be
    expressed as the quotient of a smooth manifold by an almost free
    action of a compact Lie group (see section 2). Therefore,
     we will use methods from equivariant
    topology to provide an effective K-theoretic approach to orbifolds.
    In particular, using an appropriate equivariant Chern character we
    obtain a decomposition theorem for {\sl orbifold K-theory}
    as a ring, and we show that it is additively isomorphic to orbifold
    cohomology. A nice by-product of our orbifold
    K-theory is a natural notion of orbifold Euler number for
    a general orbifold. In the context of orbifold cohomology, it only makes
    sense for $SL$-orbifolds, where the
    local group is a finite subgroup of $SL(n, C)$.

    In order to compare orbifold K-theory with orbifold cohomology,
    we make use of equivariant Bredon cohomology with coefficients in
    the representation ring functor. This equivariant theory
    is the natural target for Chern characters and with real
    coefficients is additively isomorphic to orbifold cohomology.
    It can however be defined with very general coefficients and hence
    seems to be an important technical device for the study of
    orbifolds.

    Orbifold K-theory has a different and  interesting ring
    structure. It suggests that orbifold K-theory should be interesting in its own right,
    and not just orbifold cohomology with a $\mathbb Z/2$ grading.
    This is indeed the case for several reasons. For example,
    orbifold K-theory (not orbifold cohomology) is the natural
    set-up for index theory. There is also an indication that
    orbifold K-theory behaves better algebro-topologically than
    orbifold cohomology.
    In \cite{CR1}, Chen-Ruan introduced the notion of a ``good map''
    to replace a traditional orbifold map as an
    appropriate map between
    orbifolds. We shall refer to these as \emph{orbifold morphisms};
    an equivalent notion appreared previously in \cite{MP1}.
    From their definition it is not hard to
    see that an orbifold morphism is natural with respect to orbifold
    K-theory and that it induces a ring
    homomorphism on orbifold K-theory. 
    In contrast, such a
    naturality property
    is far more subtle and unknown for orbifold cohomology.

    A key physical concept of orbifold string theory is the twisting by
    discrete torsion. In a cohomological context, it motivated the
    construction of twisted orbifold cohomology \cite{R}.
    Another important goal of this paper is to introduce twisting
    for orbifold K-theory. For twisted orbifold cohomology, the
    second author introduced the more general notion of an
    ``inner local
    system'' to twist orbifold cohomology. We do not know how to
    twist orbifold K-theory in such a generality; here, we
    restrict ourselves to the use of discrete torsion.
    We introduce {\sl twisted orbifold K-theory} using an explicit
    geometric model. In the case when the orbifold is of the form
    $X=Y/G$, where $G$ is a finite group, then our construction
    can be understood as a twisted version of equivariant K-theory,
    where the twisting is done using a fixed element $\alpha\in
    H^2(G,\mathbb S^1)$. The basic idea is to use the associated
    central extension, and to consider equivariant bundles with
    respect to this group which cover the $G$--action on $Y$. A computation
    of the associated twisted theory can be explicitly obtained
    (over the complex numbers) using ingredients
    from the classical
    theory of projective representations.  In fact,  we define a
    bigraded ring denoted the {\sl total twisted K-theory}
    which is parametrized over all discrete torsion elements.
    This construction seems to be of independent interest.

    More generally we can define a twisted orbifold K-theory
    associated to the universal orbifold cover; in this generality
    it can be computed in terms of twisted Bredon cohomology. This can
    be understood as the $E_2$--term of the twisted version of a spectral
    sequence converging to twisted orbifold K-theory, where
    in all known instances the
    higher differentials are trivial in characteristic zero
    (this is a standard observation in the case of the
    Atiyah--Hirzebruch spectral sequence). In particular, we obtain
    an additive isomorphism between our twisted orbifold K-theory and
    twisted orbifold cohomology for the case of global quotients.

    It remains to relate our construction to other forms of twisting
    for K-theory (such as that described in \cite{W}).
    For "good" orbifolds, a
    twisted version of K-theory has been studied by Marcolli-Mathai \cite{MM}
    using $C^*$ algebras.
    However, we would like to remark that
    the connection from geometric K-theory to algebraic K-theory
    of $C^*$ algebras
    breaks down for non-reduced orbifolds, where
    an orbifold bundle does not have any nonzero section in
    general.

    The second author would like to thank R. Dijkgraaf, W. Wang, E. Witten and
    E. Zaslow for
    many stimulating discussions.

\section{Orbifolds and Group Actions}

    Our basic idea to study orbifold K-theory is to apply methods
    from equivariant topology.
    In this section, we collect some basic
    properties of orbifolds and describe how they
    relate to group actions. For the reader's convenience,
    we include a definition of an
    orbifold and of orbifold vector bundles. One can find
    more details and alternative definitions in \cite{CR1}. 
The notion of an orbifold was first introduced by Satake in \cite{S},
where a different name, $V$-manifold, was used. Satake's
V-manifold corresponds to a reduced orbifold in our case.

    \vskip 0.1in
    \begin{definition}

An orbifold structure on a Hausdorff, separable
topological space $X$ is given by an
open cover ${\U}$ of $X$ satisfying the following conditions:
\begin{itemize}
\item {(2.1.1a)} Each element $U$ in ${\U}$ is uniformized, say by
$(V,G,\pi)$. Namely, $V$ is a smooth manifold and $G$ is a finite
group acting smoothly on $V$ such that $U=V/G$ with $\pi$ as
projection map. Let $ker~G$ be the subgroup of $G$ acting
trivially on $V$.
\item {(2.1.1b)} If $U^\prime\subset U$, then there is a collection of
injections $(V^\prime,G^\prime,\pi^\prime)\rightarrow (V,G,\pi)$.
Namely, the inclusion $i: U^\prime\subset U$ can be lifted to maps
$\tilde{i}: V^\prime \rightarrow V$ and an injective homomorphism
$i_{\#}: G^\prime \rightarrow G$ such that $i_{\#}$ is an
isomorphism from $ker~G^\prime$ to $ker~G$ and  $\tilde{i}$ is
$i_{\#}$-equivariant.
\item {(2.1.1c)} For any point $p\in U_1\cap U_2$, $U_1,U_2\in {\U}$,
there is a $U_3\in {\U}$ such that $p\in U_3\subset U_1\cap U_2$.
\end{itemize}.
    \end{definition}

   For any point $p\in X$, suppose that $(V,G,\pi)$ is a
   uniformizing neighborhood and $\bar{p}\in \pi^{-1}(p)$. Let $G_p$
   be the stabilizer of $G$ at $\bar{p}$. Up to conjugation, it is
   independent of the choice of $\bar{p}$ and called the {\em local
   group of $p$}.
   \begin{definition}
      We call $X$ {\em a reduced orbifold} if $G_p$ acts effectively
      for all $p\in X$.
      \end{definition}

    It is well-known that if a compact Lie group $G$ acts smoothly
    and effectively
    on a manifold $M$ with finite stabilizers (isotropy subgroups),
    then
    $M/G$ is a reduced orbifold. More generally, $X=M/G$ is an orbifold
    for any smooth Lie group action $G$ if the following conditions are
    satisfied:

(2i) For any $x\in M$ the isotropy subgroup $G_x$ is finite (this
is what we will call an \emph{almost free} $G$--action).
\medskip

(2ii) For any $x\in M$ there is a smooth slice  $S_x$ at $x$.
\medskip

(2iii) For any two points $x,y\in M$ such that $y\notin Gx$, there
are slices $S_x$ and $S_y$ such that $GS_x\cap GS_y=\emptyset$.
\bigskip

If $G$ is compact, an almost free $G$-action automatically
satisfies (2ii), (2iii). Examples arising from proper actions of
discrete groups will also appear in our work.

    Conversely, we do not know if every
    orbifold can be expressed as
    a quotient of a smooth manifold by an almost free
    action of a Lie group. We call
    such an orbifold a {\em quotient orbifold}.\footnote{Traditionally
    orbifolds of the form $M/G$, where $G$ is a finite group acting
    smoothly on a manifold $M$ have been called global quotients, we will
    adhere to this.}
    It is
    well-known that every reduced orbifold is a quotient:

\begin{proposition}
Let $X$ denote an $n$--dimensional reduced orbifold. There exists
a smooth, effective, almost free
$G=O(n)$ action on a smooth manifold $M$ such that the
quotient space $M/G$ is naturally isomorphic as an orbifold to
$X$.
\end{proposition}

In other words we can think of any reduced orbifold as the
quotient of a manifold under a certain $G=O(n)$ action. An
explicit construction for $M$ is the so-called frame bundle over
$X$, described as follows. We can choose an orbifold Riemannian
    metric on $X$ and consider the frame bundle $P_X$ of the metric.
    When $X$ is reduced, $P_X$ is a smooth manifold.
Indeed if $(V,G_x)$ is a local chart for $X$ at $x$, then
above it we have the manifold $V\times_{G_x}G$, and the quotient by the
right $G$ action gives $V/G_x$. From our point of view the above can
be interpreted as saying that the category of
reduced orbifolds can be
embedded in the category of smooth $O(n)$--manifolds with the three
properties mentioned above; a section is obtained by taking $O(n)$
orbit spaces.

It is not hard to show that if $M/G$ is a quotient orbifold,
then there exists a group extension $1\to G_0\to G\to G_{eff}\to 1$
where $G_0$ is a finite group
which is the kernel of the action (i.e. $G_0$ acts trivially on $M$)
and $G_{eff}$ acts effectively on $M$. Hence we see that we
can associate a canonical \emph{reduced} orbifold,
$X_{red}=M/G_{eff}$, to any quotient orbifold $X=M/G$.

We will assume for simplicity that our orbifolds are compact. In the
case of quotient orbifolds $M/G$ with $G$ a compact Lie group, this
is equivalent to the compactness of $M$ itself (see \cite{Bred},
page 38); a fact which we will make use of.

Next we recall the notion of the singular set of an orbifold and
its resolution.

\begin{definition}
Given an orbifold $X$, its \emph{singular set}, $\Sigma X$ consists
of the set of points $p\in X$ such that the local group $G_p\ne 1$.
The \emph{orbifold resolution} of the singular set is defined as 
$\widetilde\Sigma X = \{ (p,(g)_{G_p})~|~p\in\Sigma X,~1\ne g\in
G_p\}$, where $(g)_{G_p}$ denotes the conjugacy class of $g$ in
$G_p$.
\end{definition}

Note that we have a map $\widetilde\Sigma X
\to X$, $(p,(g))\mapsto p$.
It turns out that $\widetilde\Sigma X$ is naturally an orbifold,
which is not
necessarily connected. Later we shall see that invariants of an orbifold
$X$ can be expressed in terms of topological invariants of its orbifold
resolution.

In order to apply methods from algebraic topology in the study of orbifolds,
we recall a well known result about manifolds with smooth actions
of compact Lie groups (see \cite{Il}):
\begin{theorem}
If a compact Lie group $G$ acts on a smooth, compact manifold $M$, then
the manifold is triangulable as a finite $G$--CW complex.
\end{theorem}

Hence any such manifold will have a cellular $G$ action such that
the orbit space $X/G$ has only finitely many cells.

    We recall the notions of an orbifold vector bundle and
    of an orbifold morphism. 
    \vskip 0.1in
    \begin{definition}
     An orbifold vector bundle over an orbifold $X$ consists of
    the following data: A compatible cover ${\U}$ of $X$ such that
for any injection $i:(V^\prime,G^\prime,\pi^\prime)\rightarrow
(V,G,\pi)$, there is a continuous map $g_i:V^\prime\rightarrow
Aut(\mathbb C^k)$ giving an open embedding $V^\prime\times\mathbb C^k\rightarrow
V\times\mathbb C^k$ by $(x,v)\rightarrow (i(x),g_i(x)v)$, and for any
composition of injections $j\circ i$, we have $$ g_{j\circ
i}(x)=g_j(i(x))\circ g_i(x), \forall x\in V. \leqno (2.1.2) $$ Two
collections of maps $g^{(1)}$ and $g^{(2)}$ define isomorphic
bundles if there are maps $\delta_V:V\rightarrow Aut(\mathbb C^k)$ such
that for any injection $i: (V^\prime,G^\prime,\pi^\prime)
\rightarrow (V,G,\pi)$, we have $$
g^{(2)}_i(x)=\delta_V(i(x))\circ g^{(1)}_i(x)\circ
(\delta_{V^\prime}(x))^{-1}, \forall x\in V^\prime. \leqno (2.1.3)
$$
\end{definition}

Since $(2.1.2)$ behaves naturally under constructions of vector
spaces such as tensor product, exterior product, etc. we can
define these constructions for orbifold
vector bundles.
    Suppose that $X=M/G$ as defined previously and that
    $E\rightarrow M$
    is a $G$-equivariant bundle. Then, $E/G\rightarrow X$ is an
    orbifold vector bundle.

    Recall that if $p: X\rightarrow Y$ is a continuous map and $E$
    is a vector bundle over $Y$, then $p^*E$ is a vector bundle
    over $X$. This naturality property breaks down for orbifold
    maps. This suggests that orbifold map is a bad notion.
    One should replace it by a
    ``good'' map, as introduced by Chen-Ruan \cite{CR1}, which we prefer
to call an \emph{orbifold morphism}.\footnote{We should also point out that
an equivalent notion has been defined in \cite{MP1}, pp.~17 in the language
of groupoids; they are referred to as \emph{strong maps}.}  

    \begin{definition}
        Let $\tilde{f}:X\rightarrow X^\prime$ be a $C^\infty$ map between
orbifolds $X$ and $X^\prime$ whose underlying continuous map is
denoted by $f$. Suppose there is a compatible cover ${\U}$ of $X$,
and a collection of open subsets ${\U}^\prime$ of $X^\prime$
satisfying $(2.1.1a-c)$ and the following condition: There is a
1:1 correspondence between elements of $\U$ and ${\U}^\prime$, say
$U\leftrightarrow U^\prime$, such that $f(U)\subset U^\prime$, and
an inclusion $U_2\subset U_1$ implies an inclusion
$U_2^\prime\subset U_1^\prime$. Moreover, there is a collection of
local $C^\infty$ liftings $\{\tilde{f}_{UU^\prime}\}$ of $f$,
where $\tilde{f}_{UU^\prime}:(V,G,\pi)\rightarrow
(V^\prime,G^\prime,\pi^\prime)$ satisfies the following condition:
each  injection $i:(V_2,G_2,\pi_2)\rightarrow (V_1,G_1,\pi_1)$ is
assigned  an injection
$\lambda(i):(V_2^\prime,G_2^\prime,\pi^\prime_2)\rightarrow
(V^\prime_1,G_1^\prime,\pi_1^\prime)$ such that
$\tilde{f}_{U_1U^\prime_1}\circ
i=\lambda(i)\circ\tilde{f}_{U_2U^\prime_2}$, and for any
composition of injections $j\circ i$, the following compatibility
condition holds: $$ \lambda(j\circ i)=\lambda(j)\circ \lambda(i).
\leqno (2.2.1) $$ We call $\{\tilde{f}_{UU^\prime},\lambda\}$ a
{\it compatible system} for $\tilde{f}$.
 \end{definition}

    Given a compatible system $\xi$, one can construct a pull-back
    orbifold vector bundle in the following sense. Suppose that
    $E\rightarrow Y$ is an orbifold vector bundle. By the
    definition, $E$ can be viewed as a collection of local
    orbifold vector bundles satisfying (2.1.2), (2.1.3). One can
    pull back the local orbifold bundles. The conditions of
    a compatible system ensures that one can glue them together to
    form a global orbifold bundle denoted by $f^*_{\xi}E$. It is
    clear that there is a bundle map $\bar{f}:
    f^*_{\xi}E\rightarrow E$ covering $f$.
    \begin{definition}
      Two compatible systems
$\xi_i$ for $i=1,2,$ of a $C^\infty$ map $\tilde{f}:X\rightarrow
X^\prime$ are {\it isomorphic} if for any orbifold vector bundle
$E$ over $X^\prime$ there is an isomorphism $\psi$ between the
corresponding pull-back bundles $f^*_{\xi}E^\ast_i$ with
$\bar{f}_i:f^*_{\xi}E_i^\ast\rightarrow E$, $i=1,2$, such that
$\bar{f}_1=\bar{f}_2\circ\psi$.
    \end{definition}

    We often drop $\xi$ from the definition without any confusion.

    \begin{definition}
     An \emph{orbifold morphism} $f$ is an orbifold map together with an isomorphism
     class of compatible systems.
    \end{definition}
Later we will see that orbifold morphisms are compatible with respect to
products in K-theory, making them especially significant.

\section{Equivariant Cohomology and Orbifolds}

In this paper, we only consider quotient orbifolds $X=M/G$. We
will consider actions of both compact and discrete groups.
Throughout the paper, we will use $G$ to denote a compact Lie
group and $\Gamma$ to denote a discrete group. By the previous
discussion, a reduced orbifold is a quotient orbifold through the frame
bundle construction. Hence, our results work for a rather general
class of orbifolds. It is clear
that one can study them using equivariant techniques. We begin by
considering the equivariant cohomology of $M$. Recall that the
Borel construction associated to the $G$--action on $M$ is the
orbit space $M\times_{G}EG=(M\times EG)/G$, where $G$ is the
universal $G$--space. The $G$--equivariant cohomology of $M$ is by
definition the cohomology of the associated Borel construction. We
begin by computing this equivariant cohomology under certain
assumptions on the isotropy subgroups.
\begin{theorem}
Suppose that $X=M/G$ be a quotient orbifold. If $R$ is a ring such
that $|G_x|$ is invertible in $R$ for each $x\in X$, then there is
an algebra isomorphism $H^*_{G}(M,R)\cong H^*(X,R).$
\end{theorem}
\begin{proof}
The proof is a straightforward application of the Leray spectral
sequence of the map $\phi :M\times_{G}EG\to M/G=X$, noting that
the fibers are spaces of the form $BG_x$, which have trivial
reduced cohomology with $R$--coefficients.
\end{proof}

Consider now the so-called classifying
map
$p: M\times_{G}EG\to BG$
induced by projection onto the second factor; here $BG=EG/G$
denotes the classifying space of $G$. Assuming that $R$ is a ring
where all the $|G_x|$ are invertible, we obtain an algebra
homomorphism
$q= (\phi^*)^{-1}p^*: H^*(BG, R)\to H^*(X,R)$.
We now use this to define characteristic classes for the orbifold
$X$.

\begin{definition}
Let $X=P_X/O(n)$ be a reduced orbifold such that all its local
transformation groups are of odd order. The $i$--th
Stiefel--Whitney class $w_i(X)\in H^i(X,\mathbb Z/2)$ is defined
as $w_i(X)=q(w_i)$, where $w_1,\dots ,w_n \in H^*(BO(n),\mathbb
Z/2)$ are the usual universal Stiefel--Whitney classes.
\end{definition}
Similarly we can define Chern classes as follows

\begin{definition}
Let $X=P_X/U(n)$ be a reduced complex orbifold. 
Let $R$ denote a
ring where all the orders of the local transformation groups are
invertible. The $i$--th Chern class $c_i(X)\in H^{2i}(X,R)$ is
defined as $c_i(X)=q(c_i)$, where $c_1,\dots , c_n\in
H^*(BU(n),R)$ are the usual universal Chern classes and $q$ is the
analogue of the map previously defined.
\end{definition}

An appropriate ring $R$ can be constructed from the integers by inverting the
least common multiple of the orders of all the local transformation groups;
the rational numbers $\mathbb Q$ are of course also a suitable choice.

More generally what we see is that with integral coefficients, the
equivariant cohomology of $M$ will have interesting \emph{torsion}
classes. Unfortunately, integral computations are notoriously
difficult, especially when finite group cohomology is involved.
The mod $p$ equivariant cohomology of the frame bundle $M$ will
contain interesting information about the action; in particular
its Krull Dimension will be equal to the maximal $p$--rank of the
isotropy subgroups (see \cite{Q}). However for our geometric
applications it is convenient to use an equivariant cohomology
theory which has substantial {\sl torsion--free} information. That
is where K-theory naturally comes in, as instead of cohomology the
basic object is a representation ring.

Less well known than ordinary equivariant cohomology is
the \emph{Bredon cohomology} associated to a group action.
It is in fact the most adequate equivariant cohomology theory
available. We briefly sketch its definition for the case of
compact Lie groups (we refer to \cite{Bre}, \cite{LO}, \cite{H}
and \cite{JMO}, appendix).

Let $Or(G)$ denote the homotopy category whose objects are orbits
$G/H$, $H\subset G$ and where $Mor_{Or(G)}(G/H, G/K)$ is the set
of $G$--homotopy classes of $G$ maps between these orbits. A
\emph{coefficient system} for Bredon cohomology is a functor
$F:Or(G)^{op}\to Ab$. Now for any $G$--CW complex $M$, define
$C_*^G(M):Or(G)\to Ab_*$ by setting
$C^G_*(M)(G/H)=C_*(M^H/WH_0)$. Here $C_*(-)$ denotes the cellular
chain complex and $WH_0$ is the identity component of $NH/H$. We
now define $C^*_G(M;F)=Hom_{Or(G)}(C_*^G(M),F)$ and
$H^*_G(M;F)=H(C^*_G(M;F))$. One can see that in fact for each
$n\ge 0$, $C_G^n(M;F)$ is the direct product, taken over all
orbits $G/H\times D^n$ of $n$--cells in $M$, of the groups
$F(G/H)$. We have that $C^*_G(M;F)=Hom_{Or(G)}(C^G_*(M), F)$;
indeed this is determined on $Or(G,M)$, the full subcategory
containing all orbit types in $M$. Note that from its definition
there will be a spectral sequence (see \cite{H})
$$E_2 = Ext^*_{Or(G)}(\underline{H}_*(M), F)\Rightarrow
H^*_G(M;F)$$ where $\underline{H}_*(M)(G/H) = H_*(M^H/WH_0,
\mathbb Z)$.

In our applications we will assume that the isotropy is
all finite. Our basic example will involve the
complex representation ring functor $R(-)$; i.e.
$G/H\mapsto R(H)$. In this case the fact that $R(H)$ is a
ring implies that Bredon cohomology will have a natural
ring structure (constructed using the diagonal).

We will be interested in using the rationalized functor
$R_Q= R(-)\otimes \mathbb Q$. For $G$ finite, it was
shown in \cite{S} that $R_Q$ is an injective functor;
similarly for $\Gamma$ a discrete group and for proper
actions with finite isotropy it was shown in \cite{LO}
that $R_Q$ is injective. This result will also hold
for $G$--CW complexes with finite isotropy, where
$G$ is a compact Lie group. This follows by adapting
the methods in \cite{Sl} and described in \cite{H}.
The key technical ingredient is the surjectivity of
$R_Q(H)\to \lim_{K\in F_p(H)}~R_Q(K)$
where $H$ is any finite subgroup of $G$ and $F_p(H)$
denotes the family of all proper subgroups in $H$.
Hence we have the following basic
isomorphism:
$H^*_G(M,R_Q) \cong Hom_{Or(G)}(\underline{H}_*(M), R_Q)$.
As we shall see, in the case of a quotient orbifold $M/G$, its
orbifold cohomology is {\sl additively} isomorphic to
$H^*_G(M,R(-)\otimes\mathbb R)$.

    Suppose that $X=M/G$ is a quotient orbifold. Using
    equivariant K-theory we will show that the
    Bredon cohomology $H^*_G(M, R_Q)$ is independent of the
    presentation $M/G$ and canonically associated with
    the orbifold
    $X$ itself. A direct proof with more general coefficients
    would be of some interest.
    In the case of a reduced orbifold, we can canonically
    associate to it the Bredon cohomology of its
    frame bundle; motivated by this we
    introduce the following

\begin{definition}
Let $X$ be a reduced orbifold; we define its orbifold Bredon
cohomology with $R_Q$--coefficients as
$ H^*_{orb}(X, R_Q)=H^*_{O(n)}(P_X, R_Q)$.
\end{definition}

\begin{remark}
The approach in \cite{MP1} and \cite{MP2} using groupoids and
their classifying spaces is also an important and basic point
of view in studying orbifolds. From our perspective, if
$\mathcal{G}$ is the groupoid associated to the quotient
orbifold $X=M/G$,
then in fact $B\mathcal{G}\simeq M\times_GEG$, where 
$B\mathcal{G}$ denotes the \emph{classifying space} of the
groupoid $\mathcal{G}$. Hence
the `cohomology of the orbifold' can be identified with the
usual equivariant cohomology, and the Leray spectral sequence
of the projection onto $X=M/G$ gives rise to the basic
calculational
device for its computation. In fact the $E_2$ term of this spectral
sequence can be identified with the Bredon cohomology of the
$G$--space $M$ with coefficients in the graded functors 
$G/H\mapsto H^*(H,\mathbb Z)$, the ordinary group cohomology.
Note that if $X$ is an orbifold
which can be presented as a quotient in two different ways, say
$X\cong M/G\cong M'/G'$, then we have homotopy equivalences 
$B\mathcal{G}\simeq M\times_GEG\simeq M'\times_{G'}EG'$, hence
the equivariant cohomology is an intrinsic invariant of a quotient
orbifold. Similarly one can define the \emph{orbifold fundamental
group}, denoted $\pi_1^{orb}(X)$ as $\pi_1(B\mathcal{G})$; for a quotient orbifold
$X=M/G$ it will be isomorphic to
$\pi_1(M\times_GEG)$.
We emphasize that the `cohomology
of the orbifold $X$' is very different from the `orbifold
cohomology of $X$'. Indeed as we shall see the latter is
actually additively isomorphic to the orbifold Bredon cohomology
defined above.
\end{remark}

\section{Orbifold Bundles and Equivariant K-Theory}

We will recall (see \cite{SW}) an invariant constructed from
orbifold vector bundles, as defined in \S 2.

\begin{definition}
Given a compact orbifold $X$ we define $K_{orb}(X)$ to be the
Grothendieck group of isomorphism classes of orbifold vector
bundles on $X$.
\end{definition}

 Recall that an orbifold morphism $f:X\to Y$ is an
orbifold map together with an isomorphism
     class of compatible system.
    From the definition, one can verify that orbifold bundles over
$Y$ pull back to orbifold bundles over $X$ under an orbifold morphism. In fact
one can establish the

    \begin{proposition}
     Suppose that $(f, \xi): X\rightarrow Y$ is an
    orbifold morphism. Then, it induces a ring homomorphism
    $f^*_{\xi}: K_{orb}(Y)\rightarrow K_{orb}(X)$.
    \end{proposition}

    An important example of an orbifold morphism is the projection map $p: M\rightarrow
    M/G$, where $G$ is a compact Lie group acting
    almost freely on a manifold $M$.
    Therefore, if $E$ is an orbifold vector bundle over
    $M/G$, $p^*E$ is a
    smooth vector bundle over $M$. It is obvious that $p^*E$ is
    $G$-equivariant. Conversely, if $F$ is a $G$-equivariant bundle
    over $M$, $F/G\rightarrow M/G$ is an orbifold vector bundle
    over $X=M/G$. Therefore, we have a canonical identification
    between $K_{orb}(X)$ and $K_G(M)$.

    \begin{proposition}
     Let $X=M/G$ be a quotient orbifold. Then the projection
    map $p: M\rightarrow X$ induces an isomorphism
    $p^*: K_{orb}(X)\rightarrow K_G(M)$.
    \end{proposition}

In particular, if $X$ is a reduced orbifold, we can identify its orbifold
K-theory with the
equivariant K-theory of its frame bundle.\footnote{This has also been proposed by
J.Morava \cite{M} and it also appears implicitly in \cite{SW}}
Based on the above we can now introduce an alternative definition
of orbifold K-theory for quotient orbifolds.

\begin{definition}
Let $X=M/G$ denote a compact quotient orbifold, then
we define its orbifold K-theory as
$K^*_{orb}(X)=K^*_{G}(M)$.
\end{definition}

Note that by Bott periodicity, this invariant is $\mathbb
Z/2$--graded. It is of course possible to extend the original
definition of orbifold K-theory in the usual way; indeed if $X$ is
an orbifold, then $X\times\mathbb S^n$ is also an orbifold and
moreover the inclusion $i:X\to X\times \mathbb S^n$ is an orbifold morphism.
Let $i^*_n: K_{orb}(X\times\mathbb S^n)\to K_{orb}(X)$; then we
can define $K_{orb}^{-n}(X)= ker~i^*_n$. However the canonical
identification outlined above shows that for a quotient
orbifold this extension must agree
with the usual extension for equivariant complex K-theory to a
$\mathbb Z/2$--graded theory. Given this, we have chosen to {\sl
define} orbifold K-theory by using equivariant K-theory, as it
will enable us to make some meaningful computations. Note that if
an orbifold $X$ is presented in two different ways as a quotient
orbifold, say $M/G\cong X\cong M'/G'$, then we have shown that
$K^*_{orb}(X)\cong K^*_G(M)\cong K^*_{G'}(M')$.
Another point to make is that the homomorphism
$G\to G_{eff}$ will induce a ring map $K^*_{orb}(X_{red})
\to K^*_{orb}(X)$.

We also introduce the (K-theoretic) orbifold Euler characteristic
\footnote{This definition extends the string-theoretic orbifold
Euler characteristic which has been defined for global quotients.}

\begin{definition}
The orbifold Euler characteristic of $X$ is
\[
\chi_{orb}(X)=dim_{\mathbb Q}~K^0_{orb}(X)\otimes\mathbb Q
~-~dim_{\mathbb Q}~K^1_{orb}(X)\otimes\mathbb Q
\]
\end{definition}

It remains to show that these invariants are tractable or even well--defined.

\begin{proposition}
If $X=M/G$ is a compact quotient orbifold for a compact Lie group
$G$, then $K^*_{orb}(X)$ is a finitely generated abelian group,
and the orbifold Euler characteristic is well-defined.
\end{proposition}
\begin{proof}
We know that $M$ is a finite, almost free $G$--CW complex. It
follows from \cite{Se} that there is a spectral sequence
converging to $K_{orb}(X)=K_G(M)$, with $E_1^{p,q}$--term equal to
zero if $q$ is odd and equal to $\bigoplus_{\overline{\sigma}\in
X^{(p)}} R(G_{\sigma})$ where $X^{(p)}$ denotes the collection of
$p$--dimensional cells in $X$ and $R(G_{\sigma})$ denotes the
complex representation ring of the stabilizer of $\sigma$ in $M$.
In fact the $E_2$ term is simply the homology of a chain complex
assembled from these terms. By our hypotheses, each $G_{\sigma}$
is finite, and there are finitely many such cells, hence each term
is finitely generated as an abelian group and there are only
finitely many of them. We conclude that $E_1$ satisfies the
required finiteness conditions, and so must its subquotient
$E_{\infty}$, whence the same holds for $K^*_{orb}(X)=K^*_{G}(M)$.
\end{proof}

\begin{corollary}
$$\chi_{orb}(X) = \sum_{\overline{\sigma}\in X}~
(-1)^{dim~\overline{\sigma}}rank~R(G_{\sigma})$$ where the
$\overline\sigma$ range over all cells in $X$, with corresponding
isotropy subgroup $G_{\sigma}$.
\end{corollary}

The spectral sequence used above is in fact the equivariant
analogue of the Atiyah--Hirzebruch spectral sequence. We have
described the $E_1$--term as a chain complex assembled from the
complex representation rings of the isotropy subgroups. The
$E_2$--term coincides with the equivariant \emph{Bredon
cohomology} of $M$ with coefficients in the representation ring
functor, denoted $H^*_G(M, R(-))$. We shall see later that in fact
this spectral sequence collapses rationally at the $E_2$--term
(see \cite{LO}, page 28). As a consequence of this we have that
$H^*_{orb}(X)$, $K^*_{orb}(X)\otimes\mathbb R$ and $H^*_{G}(M,
R(-)\otimes\mathbb R)$ are all additively isomorphic and in fact
the last two invariants have the same ring structure (provided we
take the $\mathbb Z/2$--graded version of Bredon cohomology).

Computations for equivariant K-theory can be quite complicated.
Our approach will be to study the case of global quotients arising from actions of
finite and more generally discrete groups. The key computational tool will be an
equivariant Chern character, which we will define for almost free actions of compact
Lie groups. This will be used to establish the additive rational
equivalences outlined above. However it should be noted that $K_{orb}(X)$ can contain
important {\sl torsion}
classes and its rationalization is a rather crude approximation.

Let us review the special
case of a global quotient,
where the K-theoretic invariant above is more familiar.

\begin{example}
Let $G$ denote a finite group acting on a manifold $Y$ so that
$X=Y/G$ is an orbifold (a {\sl global quotient}). In this case we
know that there is an isomorphism $K_{orb}(X)\cong K_G(Y)$.
When tensored with the rationals this equivariant K-theory decomposes
as a direct sum
and we obtain the well--known formula
\[
K^*_{orb}(X)\otimes\mathbb Q \cong \bigoplus_{\{ (g)~|~g\in G\}}
K^*(Y^{\langle g\rangle}/Z_G(g))\otimes\mathbb Q
\]
where $(g)$ is the conjugacy class of $g\in G$ and $Z_G(g)$ denotes the
centralizer of $g$ in $G$. Note that this decomposition 
appears in \cite{AS} but it can be traced back (independently)
to \cite{Sl}, \cite{tD} and \cite{Ku}.
\end{example}

One of the key elements in the theory of orbifolds is the
orbifold resolution of the singular set. For the case of a global
quotient $X=Y/G$,
it can be shown (see \cite{CR}) 
that we have a homeomorphism

$$\widetilde\Sigma X\cong \bigsqcup_{\{ (g)~|~g\ne 1\} }Y^{\langle g\rangle}/Z_G(g)$$
whence in fact we see that
$K^*_{orb}(X)\cong_{\mathbb Q}
K^*(\widetilde\Sigma X\sqcup X).$
The conjugacy classes are used to index the so-called `twisted
sectors' arising
in this decomposition.
We will use this as a model for our result in the following section.

\section{A Decomposition for Orbifold K-theory}

We will now prove a decomposition for
orbifold K-theory using the methods 
developed by L\"uck and Oliver in
\cite{LO}. 
The basic technical result we will use is the construction of an
equivariant Chern character. Cohomology will be assumed $\mathbb Z/2$--graded
in the usual way.
\begin{theorem}
 Let
$X=M/G$ be a compact, quotient orbifold where $G$ is a
compact Lie group.
Then there is an equivariant Chern character which defines a
rational isomorphism of rings
\[
K^*_{orb}(X)\cong_{\mathbb Q} \prod_{\{ (C)~|~C\subset G~cyclic\}}
[H^*(M^C/Z_G(C))\otimes\mathbb Q (\zeta_{|C|})]^{W_G(C)}
\]
where $(C)$ ranges over conjugacy classes of cyclic subgroups and
$W_G(C)=N_G(C)/Z_G(C)$, a necessarily finite group.
\end{theorem}
\begin{proof}
As has been remarked, we can assume that $M$ is a finite, almost free
$G$--CW complex.
Now, as in \cite{LO} and \cite{AS}, the main idea of the proof is to construct
a natural Chern character for any $G$--space as above, and then
prove that it induces an isomorphism for orbits of the form $G/H$, where
$H\subset G$ is finite. Using induction on the number of orbit types
and a Mayer--Vietoris sequence
will complete the proof.

To begin we recall the existence (see \cite{LO}, Prop. 3.4) of a ring homomorphism
$$\psi : K^*_{N_G(C)}(M^C)\to K^*_{Z_G(C)}(M^C)\otimes R(C);$$
in this much more elementary setting it can be defined by its value on
vector bundles, namely

$$\psi ([E]) = \sum_{V\in~Irr(C)} [Hom_C(V,E)]\otimes [V]$$
for any $N_G(C)$--vector bundle $E\to M^C$. We make use of the natural maps

$$K^*_{Z_G(C)}(M^C)\otimes R(C)\to K^*_{Z_G(C)}(EG\times M^C)\otimes R(C)
\to K^*(EG\times_{Z_G(C)}M^C)\otimes R(C)$$
as well as the Chern map

$$K^*(EG\times_{Z_G(C)}M^C)\otimes R(C)\to
H^*(EG\times_{Z_G(C)}M^C;\mathbb Q)\otimes R(C)\cong
H^*(M^C/Z_G(C);\mathbb Q)\otimes R(C)$$

Note that the isomorphism above is due to the crucial fact that
all the fibers of the projection map $EG\times_{Z_G(C)}M^C\to M^C/Z_G(C)$ are
rationally acyclic, as they are classifying spaces of finite groups.
Finally we make use of the ring map $R(C)\otimes \mathbb Q\to \mathbb Q (\zeta_{|C|})$,
with kernel the ideal of elements whose characters vanish on all generators
of $C$. Putting all of this together, and using the restriction map,
we obtain a natural ring homomorphism

$$K^*_G(M)\otimes\mathbb Q \to H^*(M^C/Z_G(C),\mathbb Q (\zeta_{|C|}))^{N_G(C)/Z_G(C)}.$$
Here we have taken invariants on the right hand side, as the image naturally
lands there. Verification of the isomorphism on $G/H$ is
an elementary consequence of the isomorphism $K^*_G(G/H)\cong R(H)$ and details are
left to the reader.
\end{proof}

\begin{corollary}
Let $X=M/G$ be a compact quotient orbifold; then
there is an additive decomposition

\[
K^*_{orb}(X)\otimes\mathbb Q =K^*_G(M)\otimes \mathbb Q\cong
\bigoplus_{\{ (u)~|~u\in G\} } K^*(M^{\langle
u\rangle}/Z_{G}(u))\otimes\mathbb Q
\]
\end{corollary}

Note that the (finite)
indexing set will consist of the $G$--conjugacy
classes of elements in the isotropy subgroups--all of 
finite order. These should be considered analogues of the `twisted
sectors' which arise in the case of global quotients.

Our immediate goal is to relate
this decomposition
to specific information about the orbifold $X$.
\begin{theorem}
Let $X=M/G$ denote a compact quotient orbifold; then
there is a homeomorphism
\[
\bigsqcup_{\{ (a)~|~ a\in G\} }M^{\langle a\rangle}/Z_{G}(a) \cong
\widetilde\Sigma X\sqcup X
\]
and in particular
$K^*_{orb}(X)\cong_{\mathbb Q} K^*(\widetilde\Sigma X\sqcup X).$
\end{theorem}
\begin{proof}

We begin by considering the
situation locally. Suppose that we have a chart in $M$ of the
form $V\times_HG$, mapping onto $V/H$ in $X$, where we assume
$H\subset G$ is a finite group. Then

$$(V\times_HG)^{\langle a\rangle}=\{ H(x,u)~|~ H(x,ua)=H(x,u)\}
=\{H(x,u)~|~uau^{-1}=h\in H,~~x\in V^{\langle h\rangle}\}
$$

Let us now define an $H$ action on $\bigsqcup_{t\in H} (V^{\langle t\rangle},t)$
by $k(x,t)=(kx,ktk^{-1})$. We define a map

$$\phi : (V\times_HG)^{\langle a\rangle}\to \bigsqcup_{t\in H}(V^{\langle t\rangle},t)/H$$
by $\phi (H(x,u))=[x,uau^{-1}]$.
We check that this is well--defined; indeed if $H(x,u)=H(y,v)$ then there
exists a $k\in H$ with $(y,v)=k(x,u)$, so $y=kx, v=ku$. This means
that $vav^{-1}=kuau^{-1}k^{-1}$ and so $[y,vav^{-1}]=[kx, kuau^{-1}k^{-1}]
=[x,uau^{-1}]$ as $k\in H$. Now suppose that $z\in Z_{G}(a)$; then
$\phi (H(x,u)z)=\phi (H(x,uz))= [x, uzaz^{-1}u^{-1}]=[x,uau^{-1}]
=\phi (H(x,u))$; hence we have a well-defined map on the orbit space

$$\overline\phi : (V\times_H G)^{\langle a\rangle}/Z_{G}(a)\to
\bigsqcup_{t\in H} (V^{\langle t\rangle},t)/H.$$
This map turns out to be injective, indeed if $(x,uau^{-1})=k(y,vav^{-1})$
for some $k\in H$, then $x=ky$ and $a=u^{-1}kvav^{-1}k^{-1}u$, hence
$u^{-1}kv\in Z_{G}(a)$ and $H(x,u)(u^{-1}kv)=H(x,kv)=H(ky,kv)=H(y,v)$.
The image of $\overline\phi$ consists of the $H$--equivalence classes
of pairs $(x,t)$ where $x\in V^{\langle t\rangle}$ and $t$ is conjugate to $a$ in $G$.

Putting this together and noting that
$(V\times_HG)^{\langle a\rangle}=\emptyset$
unless $a$ is conjugate to an element in $H$, we observe that we
obtain a homeomorphism

$$\bigsqcup_{\{ (a)~|~a\in G\}} (V\times _HG)^{\langle a\rangle}/Z_{G}(a)\cong
\bigsqcup_{t\in H}(V^{\langle t\rangle},t)/H\cong \bigsqcup_{\{(t)~|~t\in H\}}
V^{\langle t\rangle}/Z_H(t).$$
To complete the proof of the theorem it suffices to observe that
by the compatibility of charts, the local homeomorphisms
on fixed-point sets can be assembled to yield the desired global
homeomorphism on $M$.
\end{proof}
\begin{remark}
The result above in fact shows that for a compact
quotient orbifold $X=M/G$, the orbifold
resolution $\tilde\Sigma X$ is isomorphic as an orbifold to 
$\bigsqcup_{\{(a)~|~a\in G\}}
M^{\langle a\rangle}/Z_G(a)$. 
\end{remark}

\begin{corollary}
Up to regrading there is an additive isomorphism
between the orbifold cohomology of $X$ and the rationalized orbifold
K-theory of $X$.
\end{corollary}
\begin{proof}
Indeed we recall that in \cite{CR}, the
orbifold cohomology $H^*_{orb}(X,\mathbb Q)$
is defined to be additively
isomorphic to $H^*(\widetilde\Sigma X\sqcup X, \mathbb Q)$. Hence via the
Chern character map we get the isomorphism above.
\end{proof}

\begin{corollary}
~~~$\chi_{orb}(X) =\chi (X\sqcup \widetilde\Sigma X)$
\end{corollary}

\begin{example}
We will now consider the case of a weighted projective space
$\mathbb C P(p,q)$ where $p$ and $q$ are assumed to be
distinct
prime numbers. Let $\mathbb S^1$ act
on the unit sphere in $\mathbb C^2$ via
$(v,w)\mapsto (z^{p}v, z^{q}w).$
The space $\mathbb CP(p,q)$ is the quotient under this
action, and it has two singular points.
$x=[1,0]$ and $y=[0,1]$. In this case the Lie group used
to present the orbifold is $SO(2)=\mathbb S^1$ and the corresponding
isotropy
subgroups are precisely $\mathbb Z/q$ and $\mathbb Z/p$.
Their fixed point sets are disjoint circles in $\mathbb S^3$, hence
the formula for the orbifold K-theory yields
$$K^*_{orb}(\mathbb CP(p,q))\cong_{\mathbb Q}
\mathbb Q (\zeta_p)\times\mathbb Q (\zeta_q)\times \Lambda (b_2)$$
where $\zeta_p$, $\zeta_q$ are the corresponding primitive roots
of unity (compare with Corollary 2.7.6 in \cite{At}). 
More explicitly we have an isomorphism
$$K^*_{orb}(X)\otimes\mathbb Q\cong
\mathbb Q [x] /(x^{p-1}+x^{p-2}+\dots + x+1)(x^{q-1}+x^{q-2}
+\dots x+1)(x^2),$$
from which we see that 
the orbifold Euler characteristic is given by 
$\chi_{orb}(\mathbb CP(p,q))= p+q$.
\end{example}
\begin{remark}
The decomposition described above is based on entirely analogous
results for proper actions of discrete groups (see \cite{LO}).
In particular this includes the case of \emph{arithmetic orbifolds},
also discussed in \cite{A} and \cite{J}.
Let $G(\mathbb R)$ denote a semisimple $\mathbb Q$--group, and $K$
a maximal compact
subgroup. Let $\Gamma\subset G(\mathbb Q$) denote an
arithmetic subgroup. Then $\Gamma$ acts on $X= G(\mathbb R)/K$,
a space diffeomorphic to Euclidean space. Moreover if $H$ is any
finite subgroup of $\Gamma$, then $X^H$ is a totally geodesic
submanifold, hence also diffeomorphic to Euclidean space.
We can make use of the Borel-Serre
completion $\overline X$ (see \cite{BS}). 
This is a contractible space with a proper
$\Gamma$--action such that the $\overline{X}^H$ are also contractible
(we are indebted to A.Borel and G. Prasad
for outlining a proof of this
\cite{BoP})
but having a compact orbit space $\Gamma\backslash {\overline X}$.
In this case we obtain the multiplicative formula
$$K^*_{\Gamma}(X)\otimes\mathbb Q\cong K^*_{\Gamma}(\overline{X})
\otimes\mathbb Q\cong
\prod_{\{ (C)~|~C\subset\Gamma~~cyclic\} }
H^*(BZ_{\Gamma}(C),\mathbb Q(\zeta_{|C|}))^{N_{\Gamma}(C)}.
$$
This allows us to express the orbifold Euler characteristic of
$\Gamma\backslash X$
in terms of group cohomology:
$$\chi_{orb}(\Gamma\backslash X) =
\sum_{\{ (\gamma)~|~\gamma\in \Gamma~of~finite~order
\}}~~~\chi (BZ_{\Gamma}(\gamma)).$$
\end{remark}
\begin{example}
Another example of some interest is that of compact, two--dimensional
hyperbolic orbifolds. They are described as quotients of the form
$\Gamma \backslash PSL_2(\mathbb R)/SO(2)$, where $\Gamma$ is a Fuchsian
subgroup. The groups $\Gamma$ can be expressed as extensions of the
form
$$1\to \Gamma'\to \Gamma \to G\to 1$$
where $\Gamma'$ is the fundamental group of a closed orientable
Riemann surface, and $G$ is a finite group (i.e. they
are virtual surface
groups). Geometrically we have an action of $G$ on a surface
$\Sigma$
with fundamental group $\Gamma'$; this action has
isolated singular points, with cyclic isotropy.
The group $\Gamma$ is $\pi_1(EG\times_G\Sigma )$, the
fundamental group of the associated Borel construction.
Assume that $G$ acts on $\Sigma$ with $n$ orbits,
having respective isotropy groups $\mathbb Z/v_1,
\dots , \mathbb Z/v_n$ and with quotient
a surface $W$ of genus equal to $g$.
The formula then yields (compare with the
description in \cite{MM}, pg. 563)
\[
K^*_{orb}(W)\otimes\mathbb Q
\cong
\tilde{R}(\mathbb Z/v_1)\otimes\mathbb Q\times\dots
\times\tilde{R}(\mathbb Z/v_n)\otimes\mathbb Q\times
K^*(W)\otimes\mathbb Q
\]
In this expression $\tilde{R}$ denotes the reduced
representation ring, which arises because the trivial
cyclic subgroup only appears once.
 From this we see that
\[
dim_{\mathbb Q}~K^0_{orb}(W)\otimes\mathbb Q =
\sum_{i=1}^n~(v_i-1) + 2,
\,\,\,\,\,\,\,\,dim_{\mathbb Q}~K^1_{orb}(W)\otimes\mathbb Q = 2g.
\]
and so we have that
$\chi_{orb}(W)~=~ \sum_{i=1}^n(v_i-1) + \chi (W)$.
\end{example}
\begin{remark}
It has been brought to our attention that the decomposition
formula is  
analogous to a decomposition of equivariant algebraic K-theory
which appears in work of Vezzosi--Vistoli \cite{VV}
and B.Toen (see \cite{T}, page 29) in the context of
algebraic Deligne--Mumford stacks. A detailed comparison would seem
worthwhile.
\end{remark}

\begin{remark}
It should also be 
observed that the decomposition above could equally well have been
stated in terms of the computation
of Bredon cohomology mentioned previously, namely $H^*_G(M,R_Q)\cong
Hom_{Or(G)}(\underline{H}_*(M),R_Q)$ and the collapse at $E_2$ of the
rationalized Atiyah--Hirzebruch spectral sequence:
$K^*_{orb}(X)\otimes\mathbb Q \cong H^*_{G}(M, R_Q)$.
It had been shown previously that a Chern character
with expected naturality
properties inducing such an isomorphism cannot exist;
in particular
\cite{H} contains an example where such
an isomorphism is impossible. However the
example is for a circle action with stationary points, our result
\footnote{I. Moerdijk has informed us that in unpublished
work (1996), he and J.Svensson obtained essentially the same Chern character 
construction as that appearing in this paper.} 
shows that \emph{almost
free} actions of compact Lie groups do indeed give rise to appropriate equivariant
Chern characters. 
A different equivariant Chern character for abelian Lie group
actions was defined in \cite{BBM}, using a $\mathbb Z/2$--indexed de Rham cohomology
(called \emph{delocalized equivariant cohomology}). Presumably it must agree with
our decomposition in the case of almost free actions.
E.Getzler has pointed out an alternative approach involving
cyclic cohomology (see \cite{BG}).
\end{remark}

\begin{remark}
If $X=M/G$ is a quotient orbifold, then the K-theory of $M\times_GEG$ and
the orbifold K-theory are related by the Atiyah--Segal Completion Theorem
in \cite{AS2}. Considering the equivariant K-theory $K^*_G(M)$ as a module 
over $R(G)$, it states that $K^*(M\times_GEG)\cong K^*_G(M)~\widehat{}$, where
the completion is taken at the augmentation ideal $I\subset R(G)$. 
\end{remark}

\section{Projective Representations, Twisted Group Algebras and Extensions}

We will now extend many of the constructions and concepts used
previously to an appropriately \emph{twisted} setting. This twisting
occurs naturally in the framework of mathematical physics and leads
to an interesting notion of `totally twisted' equivariant notions.
In this section we will always assume that we are dealing with finite
groups, unless stated otherwise. Most of the background results which
we list appear in \cite{Ka}, Chapter III.

\begin{definition}
let $V$ denote a finite dimensional complex vector space.
A mapping $\rho: G\to GL(V)$ is called a projective representation of $G$
if there exists a function $\alpha : G\times G\to \mathbb C^*$ such
that $\rho (x)\rho (y) = \alpha (x,y)\rho (xy)$ for all $x,y\in G$
and $\rho (1) = Id_V$.
\end{definition}

Note that $\alpha$ defines a $\mathbb C^*$--valued cocycle on $G$,
i.e. $\alpha \in Z^2(G,\mathbb C^*)$. Also there is a one-to-one
correspondence between projective representations of
$G$ as above and homomorphisms from $G$ to $PGL(V)$. We will be
interested in the notion of \emph{linear equivalence} of projective
representations.

\begin{definition}
Two projective representations
$\rho_1 : G\to GL(V_1)$
and $\rho_2: G\to GL(V_2)$ are said to be
linearly equivalent if there exists
a vector space isomorphism $f:V_1\to V_2$ such that
$\rho_2(g)=f\rho_1(g)f^{-1}$ for all $g\in G$.
\end{definition}

If $\alpha$ is the cocycle attached to $\rho$, we say that
$\rho$ is an $\alpha$--representation on the space $V$.
We list a couple of basic results

\begin{lemma}
Let $\rho_i$, $i=1,2$ be an $\alpha_i$--representation on the
space $V_i$. If $\rho_1$ is linearly equivalent $\rho_2$, then
$\alpha_1$ is equal to $\alpha_2$.
\end{lemma}

It is easy to see that given a fixed cocycle $\alpha$, we can take
the direct sum of any two
$\alpha$--representations. Hence we can introduce

\begin{definition}
We define $M_{\alpha}(G)$ as the monoid of linear isomorphism classes
of $\alpha$--representations of $G$. Its associated Grothendieck
group will be denoted $R_{\alpha}(G)$.
\end{definition}

In order to use these constructions we need to introduce the notion
of a \emph{twisted group algebra}. If $\alpha :G\times G\to \mathbb C^*$
is a cocycle, we denote by $\mathbb C^{\alpha}G$ the vector space
over $\mathbb C$ with basis $\{\overline{g}~|g\in G\}$ with product
$\overline{x}\cdot \overline{y} = \alpha (x,y) \overline{xy}$
extended distributively.
One can check that $\mathbb C^{\alpha}G$ is a $\mathbb C$--algebra with
$\overline{1}$ as the identity element. This algebra is called the
\emph{$\alpha$--twisted group algebra} of $G$ over $\mathbb C$. Note that
if $\alpha (x,y)=1$ for all $x,y\in G$,
then $\mathbb C^{\alpha}G=\mathbb CG$.

\begin{definition}
If $\alpha$ and $\beta$ are cocycles, then $\mathbb C^{\alpha} G$ and
$\mathbb C^{\beta}G$ are equivalent if there exists a $\mathbb C$
algebra isomorphism
$\psi : \mathbb C^{\alpha}G\to \mathbb C^{\beta}G$
and a mapping $t:G\to \mathbb C^*$ such that
$\psi (\overline{g}) = t(g) \tilde{g}$ for all $g\in G$, where
$\{\overline{g}\}$ and $\{\tilde{g}\}$ are bases for
the two twisted algebras.
\end{definition}

We have a basic result which classifies twisted group algebras.
\begin{theorem}
We have an isomorphism between twisted group algebras,
$\mathbb C^{\alpha}G\simeq \mathbb C^{\beta}G$, if and only
if $\alpha$ is cohomologous to $\beta$; hence if
$\alpha$ is a coboundary, $\mathbb C^{\alpha}G\simeq \mathbb CG$.
Indeed, $\alpha\mapsto \mathbb C^{\alpha}G$ induces a bijective
correspondence
between $H^2(G,\mathbb C^*)$ and the set of equivalence classes of
twisted group algebras of $G$ over $\mathbb C$.
\end{theorem}

Next we recall how these twisted algebras play a role in determining
$R_{\alpha}(G)$.

\begin{theorem}
There is a bijective correspondence between $\alpha$--representations
of $G$ and $\mathbb C^{\alpha}G$--modules. This correspondence preserves
sums and bijectively maps linearly equivalent (respectively irreducible,
completely reducible) representations into isomorphic (respectively
irreducible, completely reducible) modules.
\end{theorem}

\begin{definition}
Let $\alpha\in Z^2(G,\mathbb C^*)$.
An element $g\in G$ is said to be \emph{$\alpha$--regular} if
$\alpha (g,x)=\alpha (x,g)$
for all $x\in Z_G(g)$.
\end{definition}
Note that the identity element is $\alpha$--regular for all $\alpha$.
Also one can see that $g$ is $\alpha$--regular if and only if
$\overline{g}\cdot \overline{x}=\overline{x}\cdot\overline{g}$ for
all $x\in Z_G(g)$.

If an element $g\in G$ is $\alpha$--regular, then any conjugate of $g$
is also $\alpha$--regular, hence we can speak of $\alpha$--regular
conjugacy classes in $G$. For technical purposes we also want to
introduce the notion of a `standard' cocyle; it will be a cocycle
$\alpha$ with values in $\mathbb C^*$ such that  (1) $\alpha (x,x^{-1})=1$
for all $x\in G$ and (2) $\alpha (x,g)\alpha (xg, x^{-1})=1$ for all
$\alpha$--regular $g\in G$ and all $x\in G$. Expressed otherwise, this
simply means that $\alpha$ is standard if and only if for all
$x\in G$ and for all $\alpha$--regular elements $g\in G$, we have
$\overline{x}^{-1} = \overline{x^{-1}}$ and
$\overline{x}\,\overline{g}\,\overline{x}^{-1} =
\overline{xgx^{-1}}$. It can be shown that in fact any cohomology
class $c\in H^2(G,\mathbb C^*)$ can be represented by a standard
cocycle, hence we will make this assumption from now on.

The next result is basic:
\begin{theorem}
If $r_{\alpha}$ is equal to the number of non--isomorphic irreducible
$\mathbb C^{\alpha}G$--modules, then this number is equal to
the number of distinct $\alpha$--regular conjugacy classes of $G$.
In particular $R_{\alpha}(G)$ is a free abelian group of rank
equal to $r_{\alpha}$.
\end{theorem}

In what follows we will be using cohomology classes in $H^2(G,\mathbb S^1)$,
where the $G$--action on the coefficients is assumed to be trivial.
Note that $H^2(G,\mathbb S^1)\cong H^2(G,\mathbb C^*)\cong H^2(G,
\mathbb Q/\mathbb Z)\cong H^3(G,\mathbb Z)$. We will always use standard
cocycles to represent any given cohomology class.

An element $\alpha\in H^2(G,\mathbb S^1)$ corresponds to an equivalence
class of group extensions

$$1\to \mathbb S^1\to \widetilde{G}_{\alpha}\to G\to 1$$
The group $\widetilde{G}_{\alpha}$ can be given the structure of a compact Lie
group, where  $\mathbb S^1\to \widetilde{G}_{\alpha}$ is the inclusion of a closed
subgroup.
The elements in the extension group can be
represented by pairs $\{ (g,a)~|~g\in G, a\in\mathbb S^1\}$ with the
product $(g_1,a_1)(g_2,a_2)=(g_1g_2,\alpha (g_1,g_2)a_1a_2)$.

Consider the case when $z\in Z_G(g)$; 
then we can compute the following commutator of lifts:
$$(z,1)(g,1)[(g,1)(z,1)]^{-1}=(zg,\alpha (z,g))(z^{-1}g^{-1},
\alpha (g,z)^{-1})$$
$$
=(1,\alpha (zg, (zg)^{-1})\alpha (z,g)\alpha (g,z))
=(1,\alpha (z,g)\alpha (g,z)^{-1}).$$
This computation
is independent of the choice of lifts. It can be seen that this defines
a character $L_g^{\alpha}$ for the centralizer $Z_G(g)$, via the correspondence
$z\mapsto \alpha(z,g) \alpha (g,z)^{-1}$. Note that this character is trivial
if and only if $g$ is $\alpha$--regular.

There is a 1-1 correspondence
between isomorphism classes of representations of 
$\widetilde{G}_{\alpha}$ which restrict to scalar multiplication on the central
$\mathbb S^1$ and isomorphism classes of $\alpha$-representations of $G$.
If $\psi:\widetilde{G}_{\alpha}\to GL(V)$ is such a representation then we
define an associated $\alpha$-representation via
$\rho (g) =\psi (g,1)$. Note that 
$\rho (gh)=\psi (gh,1)
=\alpha (g,h)^{-1}\psi (gh, \alpha (g,h))=\alpha (g,h)^{-1}
\psi ((g,1)(h,1))
=\alpha (g,h)^{-1}\rho (g)\rho (h)$.
Conversely, given
$\rho :G\to GL(V)$, we simply define 
$\psi (g,a)=a\rho (g)$; note that
$$\psi ((g,a)(h,b))=\psi (gh, \alpha (g,h)ab)
=ab\rho(g)\rho(h)=a\rho (g)b\rho (h)=
\psi (g,a)\psi (h,b).$$
Hence we can identify $R_{\alpha}(G)$ as the subgroup of 
$R(\widetilde{G}_{\alpha})$ generated by representations that
restrict to scalar multiplication on the central $\mathbb S^1$.

In the next section we will need an explicit understanding of
the action of $\widetilde{Z_G(g)}_{\alpha}$ on $R_{res(\alpha)}
(\langle g\rangle)$, where $res(\alpha)$ is the restriction of the
cocycle to the subgroup $\langle g\rangle$ (this restriction is cohomologous
to zero).
It is easiest to describe using the formulation above. Given a
representation $\phi$ for $\widetilde{\langle g\rangle}_{\alpha}$, 
an element $(z,a)\in \widetilde{Z_G(g)}_\alpha$, and 
$(x,b)\in \widetilde{\langle g\rangle}_{\alpha}$, we define
$(z,a)\phi (x,b)=\phi ((z,a)(x,b)(z,a)^{-1})$.
Notice that this value is precisely $L_x^{\alpha}(z)\phi(x,b)$;
this is independent of the choice of lifting and defines an action
of $Z_G(g)$.
For $x,y\in \langle g\rangle$ we have $L_x^\alpha (z)L_y^\alpha (z) = L_{xy}^\alpha (z)$,
in particular if $g^n=1$ we have $[L_g^\alpha (z)]^n=1$. The correspondence
$x\mapsto L_x^\alpha (z)$ defines a character $L^\alpha (z)$
for $\langle g\rangle$, whence
the action is best described as sending an $\alpha$-representation
$\rho$ to $L^{\alpha}(z)\rho$. Note that the evaluation 
$\phi\mapsto Tr(\phi (g,1))$ defines a 
$\mathbb C Z_G(g)$--homomorphism 
$u:R_{res(\alpha)}(\langle g\rangle)\otimes\mathbb C\to L_g^{\alpha}$.

\section{Twisted Equivariant K-theory}

We are now ready to define a twisted version of equivariant K-theory
for global quotients. We assume as before that $G$ is a finite group.
Now suppose we are given a class $\alpha\in H^2(G,\mathbb S^1)$ and
the compact Lie group extension which represents it,
$1\to\mathbb S^1\to\widetilde{G}_{\alpha}\to G\to 1$. Now let $X$ denote a finite
$G$--CW complex.

\begin{definition}
An $\alpha$--twisted $G$--vector bundle on $X$ is a
complex vector bundle $E\to X$ such that
$\mathbb S^1$ acts on the fibers through complex multiplication, 
so that the action extends to an action of $\tilde{G}_{\alpha}$
on $E$ which covers the given $G$--action on $X$.
\end{definition}

In fact $E\to X$ is a $\widetilde{G}_{\alpha}$--vector bundle, where
the action on the base is via the projection onto $G$ and the
given $G$--action. Note that if we divide out by the action
of $\mathbb S^1$, we obtain a \emph{projective} bundle over
$X$. These twisted bundles can be added, hence forming a monoid.

\begin{definition}
The $\alpha$--twisted $G$--equivariant K-theory of $X$, denoted
by $^\alpha K_G(X)$, is defined as the Grothendieck group of isomorphism classes
of $\alpha$-twisted $G$-bundles over $X$.
\end{definition}

As with $\alpha$-representations, we can describe this twisted
group as the subgroup of $K_{\widetilde{G}_{\alpha}}(X)$ generated
by isomorphism classes of bundles that restrict to multiplication
by scalars on the central $\mathbb S^1$. As the $\mathbb S^1$--action on
$X$ is trivial, we have a natural isomorphism
$K_{\mathbb S^1}(X)\cong K(X)\otimes R(\mathbb S^1)$. Composing the
restriction with the map $K(X)\otimes R(\mathbb S^1)
\to R(\mathbb S^1)$ we obtain a homomorphism
$K_{\widetilde{G}_{\alpha}}(X)\to R(\mathbb S^1)$; we can define
$^\alpha K_G(X)$ as the inverse image of the subgroup generated by
the representations defined by scalar multiplication.

Just as in non--twisted equivariant K-theory, this definition extends
to a $\mathbb Z/2$--graded theory. In fact we can define
$^\alpha K_G^0(X)=^\alpha K_G(X)$ and
$^\alpha K_G^1(X)= ker~[^\alpha K_G(\mathbb S^1\times X)
\to ^\alpha K_G(X)]$. We can also extend the 
description given above to express 
$^\alpha K^*_G(X)$ 
as a subgroup of $K^*_{\widetilde{G}_{\alpha}}(X).$

We begin by considering the case $\alpha=0$; this corresponds to the
split extension $G\times \mathbb S^1$. Any ordinary $G$--vector bundle
can be made into a $G\times\mathbb S^1$--bundle via scalar multiplication
on the fibers; conversely a $G\times\mathbb S^1$--bundle restricts to
an ordinary $G$--bundle. Hence we readily see that
$^{\alpha} K^*_G(X)=K^*_G(X)$.

Now we consider the case when $X$ is a trivial $G$--space.

\begin{lemma}
Let $X$ denote a $CW$-complex with a trivial $G$--action; then there is
a natural isomorphism
$K(X)\otimes R_{\alpha}(G)\to ^\alpha K_G(X)$.
\end{lemma}
\begin{proof}
This result is the analogue of the untwisted version (see \cite{Se}, page 133).
The natural map $R(\widetilde{G}_{\alpha})\to K_{\widetilde{G}_{\alpha}}(X)$
can be combined with the map $K(X)\to K_{\widetilde{G}_{\alpha}}(X)$ 
(which gives any vector bundle the trivial $G$--action) to yield a natural
isomorphism $K(X)\otimes R(\widetilde{G}_{\alpha})\to K_{\widetilde{G}_{\alpha}}(X)$
which covers the restriction to the $\mathbb S^1$--action; the result
follows from looking at inverse images of the subgroup generated by the scalar
representation.
\end{proof}

The inverse of the map above is given by 
$$[E]\mapsto \bigoplus_{\{[M]\in Irr(\widetilde{G}_{\alpha})\}} 
[Hom_{\widetilde{G}_{\alpha}}(M,E)]\otimes[M].$$
Note that only the $M$ which restrict to scalar multiplication on $\mathbb S^1$
are relevant--these are precisely the irreducible $\alpha$--representations.

Let $X$ be a $G$--space and $Y$ a $G'$--space, and let h:$G\to G'$ denote
a group homomorphism. If $f:X\to Y$ is a continuous map equivariant with
respect to this homomorphism we obtain a map
$^\alpha f^*: ^\alpha K_{G'}(Y)\to ^{h^*(\alpha)} K_G(X)$,
where $h^*:H^2(G',\mathbb S^1)\to H^2(G,\mathbb S^1)$ is the map induced
by $h$ in cohomology.
Let $H\subset G$ be a subgroup; the inclusion defines a restriction map
$H^2(G,\mathbb S^1)\to H^2(H,\mathbb S^1)$. In fact if $\widetilde{G}_{\alpha}$
is the group extension defined by $\alpha\in H^2(G,\mathbb S^1)$; then
$res^G_H(\alpha)$ defines the `restricted' group extension over $H$,
denoted $\widetilde{H}_{res(\alpha)}$; we have a restriction map
$^\alpha K_G(X)\to ^{res(\alpha)} K_H(X)$. 
 
Now consider the case of 
an orbit $G/H$; then we have
$^\alpha K_G(G/H)=R_{res^G_H(\alpha)}(H)$. Indeed we can identify
$K_{\widetilde{G}_{\alpha}}(G/H) = 
K_{\widetilde{G}_{\alpha}}(\widetilde{G}_{\alpha}/\widetilde{H}_{\alpha})
\cong R(\widetilde{H}_{\alpha})$, and by restricting to the representations
that induce scalar multiplication on $\mathbb S^1$ we obtain the result.

We are now ready to state a basic decomposition theorem for our twisted
version of equivariant K-theory

\begin{theorem}
Let $G$ denote a finite group and $X$ a finite $G$--CW complex. For
any $\alpha\in H^2(G,\mathbb S^1)$ we have a decomposition

$$^\alpha K^*_G(X)\otimes\mathbb C\cong \bigoplus_{\{(g)~|~g\in G\}}
(K^*(X^{\langle g\rangle})\otimes L_g^{\alpha})^{Z_G(g)}.$$
\end{theorem}
\begin{proof}
Fix the class $\alpha\in H^2(G,\mathbb S^1)$; for any subgroup
$H\subset G$,
we can associate
$H\mapsto R_{res(\alpha)}(H)$.
Note the special case when $H=\langle g\rangle$, a cyclic subgroup. As
$H^2(\langle g\rangle,\mathbb S^1)=0$, $R_{res(\alpha)}(\langle
g\rangle)$ is 
isomorphic to $R(\langle g\rangle)$.

Now consider $E\to X$, an $\alpha$--twisted bundle over $X$; it
restricts to an $res(\alpha)$--twisted bundle over
$X^{\langle g\rangle}$.
Recall that we have an isomorphism
$^{res(\alpha)} K^*_{\langle g\rangle}(X^{\langle g\rangle})\cong
K^*(X^{\langle g\rangle})\otimes R_{res(\alpha)}(\langle g \rangle)$.
Let $u:R_{res(\alpha)}(\langle g\rangle)\to L_g^{\alpha}$ denote the 
$\mathbb CZ_G(g)$--map $\chi\mapsto \chi (g)$ described previously, where the
centralizer acts on the projective representations as described above.
Then the composition

$$^\alpha K^*_G(X)\otimes\mathbb C\to 
^{res(\alpha)} K^*_{\langle g\rangle}(X^{\langle g\rangle})\otimes\mathbb C\to
K^*(X^{\langle g\rangle})\otimes R_{res(\alpha)}(\langle
g\rangle)\otimes \mathbb C
\to K^*(X^{\langle g\rangle})\otimes L_g^{\alpha}$$
has image lying in the invariants under the $Z_G(g)$--action.
Hence we can put these together to yield a map

$$^\alpha K^*_G(X)\otimes\mathbb C\to
\bigoplus_{(g)}
(K^*(X^{\langle g\rangle})\otimes L_g^{\alpha})^{Z_G(g)}.$$
One checks that this induces an isomorphism on orbits
$G/H$; the desired isomorphism follows from using induction
on the number of $G$--cells in $X$ and a Mayer-Vietoris
argument (as in \cite{AS}).
\end{proof}
 From the definition of twisted orbifold cohomology appearing
in \cite{Ru}, we conclude

\begin{corollary}
After regrading there is an additive isomorphism between the
$\alpha$--twisted
orbifold cohomology of a global quotient $M/G$, and the
$\alpha$--twisted equivariant K-theory $^\alpha K_G^*(M)
\otimes\mathbb C$.
\end{corollary}

Note that in the case when $X$ is a point, we  are saying that
$R_{\alpha}(G)\otimes\mathbb C$ has rank equal to the number
of conjugacy classes of elements in $G$ such that the associated
character $L_g^{\alpha}$ is trivial. This of course agrees with
the number of $\alpha$--regular conjugacy classes, as indeed
$^\alpha K_G(*)=R_{\alpha}(G)$.

The reader may have noticed that our twisted equivariant K-theory
does not have a product structure. Moreover it depends on a choice
of a particular cohomology class in $H^2(G,\mathbb S^1)$. Our next
goal is to relate the different twisted versions by using a product
structure inherited from the additive structure of group extensions.

Suppose we are given $\alpha, \beta$ in $H^2(G,\mathbb S^1)$,
represented by central extensions
$1\to\mathbb S^1\to\tilde{G}_1\to G\to 1$ and
$1\to\mathbb S^1\to\tilde{G}_2\to G\to 1$.
These give rise to a central extension of the form
$$1\to\mathbb S^1\times\mathbb S^1\to \tilde{G}_1\times\tilde{G}_2
\to G\times G\to 1.$$
Now we make use of the diagonal embedding $\Delta:G\to G\times G$ and
the product map $\mu:\mathbb S^1\times\mathbb S^1\to \mathbb S^1$
to obtain a central extension
$$1\to\mu (\mathbb S^1\times\mathbb S^1)\to\tilde{G}\to \Delta (G)\to 1.$$
This operation corresponds to the sum of cohomology classes, i.e.
the extension above represents $\alpha + \beta$.
Note that $ker~\mu = \{ (z, z^{-1})\}\subset\mathbb S^1\times\mathbb S^1$.

Now consider an $\alpha$--twisted bundle $E\to X$ and a $\beta$--twisted
bundle $F\to X$. Consider the tensor product bundle $E\otimes F\to X$.
Clearly it will have a $\tilde{G}_1\times\tilde{G}_2$ action on it,
which we can restrict to the inverse image of $\Delta (G)$. Now note
that $ker~\mu$ acts trivially on $E\otimes F$, hence we obtain
a $\tilde{G}$ action on $E\otimes F$, covering the $G$--action on $X$.
This is an $\alpha + \beta$--twisted bundle over $X$. Hence we have
defined a product
$$^\alpha K_G(X)\otimes^\beta K_G(X)\to ^{\alpha +\beta}K_G(X)$$
which prompts us to introduce the following definition.
\begin{definition}
The \emph{total twisted equivariant K-theory} of a $G$--CW complex
$X$ is defined as
$$TK^*_G(X) = \bigoplus_{\alpha\in H^2(G,\mathbb S^1)}~^\alpha K^*_G(X)$$
\end{definition}

Using the product above, we deduce that $TK^*_G(X)$ is a bigraded algebra,
as well as a module over $K^*_G(X)$. Note that the indexing set is finite,
and that in particular given any homogeneous element a sufficiently high
power of it will land in $^0K^*_G(X)$.

We obtain a purely algebraic construction from the above when $X$ is
a point. Namely we obtain the \emph{total twisted representation ring}
of $G$, defined as
$$TR(G)=\bigoplus_{\alpha\in H^2(G, \mathbb S^1)}~ R_{\alpha}(G),$$
endowed with
the graded algebra structure defined above. Note that if a cohomology
class is represented by a cocycle $\mu$, then its negative is
represented by $\mu^{-1}$.
Hence we see that $\rho\mapsto\rho^*$
defines an isomorphism between
$R_{\alpha}(G)$ and $R_{-\alpha}(G)$. Indeed, using vector bundles instead
we can easily extend this to show that $^\alpha K_G(X)$ is isomorphic
to $^{-\alpha}K_G(X)$.
We now provide some examples to illustrate the properties of
this construction.
\begin{remark}
It is apparent that the constructions introduced in this section can
be extended to the case of a proper action on $X$ of a discrete group
$\Gamma$. The group extensions and vector bundles used for the
finite group case have natural analogues, and so we can define
$^\alpha K^*_{\Gamma}(X)$ for $\alpha\in H^2(\Gamma, \mathbb S^1)$.
We will make use of this
in the next section.
\end{remark}

\begin{example}
Consider the group $G=\mathbb Z/2\times\mathbb Z/2$; then
$H^2(G,\mathbb S^1)=\mathbb Z/2$ (as can be seen from the Kunneth formula).
If $a, b$ are generators for $G$, we have a projective representation
$\mu: G\to PGL_2(\mathbb C)$
given by

$$a\mapsto\left(\begin{array}{cc}
               0&1\cr
           -1&0\end{array}\right), ~~~~b\mapsto\left(\begin{array}{cc}
           -1&0\cr
           0&1\end{array}\right)$$
Note that this gives rise to an extension $\tilde{G}\to GL_2(\mathbb C)$.
Restricted to $\mathbb Z/2\subset\mathbb S^1$, we get an extension of the form
$1\to \mathbb Z/2\to \tilde{D}\to\mathbb Z/2\times\mathbb Z/2\to 1$; however
this is precisely the embedding of the dihedral group in $GL_2(\mathbb C)$.
Hence the extension $\tilde{G}$ must also be non--split, and so represents
the non--trivial element $\alpha$ in $H^2(G,\mathbb S^1)$.
One can easily verify that there is only one conjugacy class of $\alpha$
regular elements in $G$, the trivial element. The representation
$\mu$ is clearly irreducible, hence up to isomorphism is the unique
irreducible $\alpha$--twisted representation of $G$. In particular,
$R_{\alpha}(G)\cong \mathbb Z \langle\mu\rangle$.
Computing the different products, we obtain
$$TR(G) = \mathbb Z[x_0,y_0, \mu_1]/(x^2=1, y^2=1,
                    x\mu =\mu, y\mu =\mu, \mu^2=1+x+y+xy) $$
The ring is $\mathbb Z/2$--graded, as we have indicated with the indices.
It has total rank equal to five. More generally, if $M$ is any compact
manifold with a $G=\mathbb Z/2\times\mathbb Z/2$--action, then one can
easily verify that
$\chi (TK^*_G(M)\otimes\mathbb C) = 6~\chi (M/G)~-~\chi (M)$.
\end{example}

\begin{example}
Let $G=\Sigma_n$, the symmetric group on $n$ symbols. Assume that
$n\ge 4$; it is well--known that in this range $H^2(G,\mathbb S^1)=
\mathbb Z/2$; denote the non--trivial class by $\alpha$. Using the
decomposition formula, one can calculate (see \cite{U})
$^\alpha K^*_{\Sigma_n}(M^n)$, where the group acts on the
$n$--fold product of a manifold $M$ by permutation of coordinates
(the quotient orbifold is the symmetric product). From this
one can recover a corrected version
of a formula which appears in \cite{Di} for twisted
symmetric products. This correction was first observed and
corrected by 
W.Wang in \cite{Wa}.
This can be expressed as follows:

$$
\sum~q^n\chi (^\alpha K^*_{\Sigma_n}(M^n)\otimes\mathbb C)
=$$
$$\prod_{n>0} (1-q^{2n-1})^{-\chi (M)}
+ \prod_{n>0} (1+q^{2n-1})^{\chi (M)}[1+\frac{1}{2}
\prod_{n>0}(1+q^{2n})^{\chi (M)}
- \frac{1}{2}\prod_{n>0}(1-q^{2n})^{\chi (M)}]
$$
\end{example}

\section{Twisted Orbifold K-theory and Twisted Bredon Cohomology }

In this section, we shift back to the orbifold point of view. In the
case of twisted orbifold cohomology, one can construct twisting
using any inner local system \cite{R}. We do not know how to twist
orbifold K-theory in such generality. Here, we have the more
limited goal of twisting orbifold K-theory using discrete torsion.

Recall that a discrete torsion $\alpha$ of an orbifold $X$ is
defined as a class $\alpha\in
H^2(\pi^{orb}_1(X), S^1)$. Here, the
orbifold
fundamental group $\pi^{orb}_1(X)$ is the group of deck
translations of the orbifold universal cover of $Y\to X$.


For example, if $X=Z/G$ is a global
quotient, the universal cover $Y$ of $Z$
 is the orbifold universal cover of $X$. In fact, if
 $Z\times_GEG$ is the Borel construction for $Z$, then we have
 a fibration sequence $Z\to Z\times_GEG\to BG$ which gives rise to
 the group extension $1\to \pi_1(Z)\to \pi_1^{orb}(X)
 \to G\to 1$; here we are
 identifying $\pi_1^{orb}(X)$ with $\pi_1(Z\times_GEG)$.
 Note that a class
 $\alpha\in H^2(G, S^1)$ induces a class
 $f^*(\alpha)$ in $H^2(\pi^{orb}_1(X), S^1)$.

 Now suppose that $X=M/G$ is a quotient manifold
 for a compact Lie group $G$ and $p: Y\rightarrow X$ is the
 orbifold universal cover. Note that $p$ is an orbifold morphism. 
 The same
 argument used in pulling back orbifold bundles
 implies that we can pull back the
 orbifold principal bundle $M\rightarrow X$ to obtain an orbifold
 principal $G$-bundle $\tilde{M}\rightarrow Y$. Furthermore,
 $\tilde{M}$ is smooth and has a free
 left $\pi^{orb}_1(X)$--action, as
well as a right $G$--action. These can be combined to yield a
left $\pi=\pi^{orb}_1(X)\times G$--action. Note that we have
$$K^*_{\pi}(\tilde{M})\cong
K^*_{G}(\tilde{M}/\pi^{orb}_1(X))=K^*_{orb}(X).$$

 Consider  a group $\pi$ of the form
$\Gamma\times G$, where $\Gamma$ is a discrete group and $G$ is a
compact Lie group. Now let $Z$ denote a proper $\pi$-complex such
that the orbit space $Z/\pi$ is a compact orbifold.
We now fix a cohomology class $\alpha\in H^2(\Gamma, \mathbb S^1)$,
corresponding to a central extension $\Gamma_{\alpha}$. From this
we obtain an extension $\tilde{\pi}_{\alpha}=\tilde{\Gamma}_{\alpha}
\times G$. We can
define the $\alpha$--twisted $\pi$--equivariant K-theory
of $Z$,
denoted $^\alpha K^*_{\pi}(Z)$ in a manner analogous to
what we did before. Namely we consider $\tilde{\pi}_{\alpha}$
bundles covering the $\pi$ action on $Z$, such that the central
circle acts by scalar multiplication on the fibers. Based
on this we can
introduce\footnote{Alternatively we could have used an equivariant
version of orbifold bundles and introduced the twisting geometrically.
This works for general orbifolds but we will not elaborate
on this here.}

\begin{definition}
Let $X=M/G$ denote a compact quotient orbifold where $G$
is a compact Lie
group, and let $Y\to X$ denote its orbifold universal cover,
with deck transformation group $\Gamma= \pi_1^{orb}(X)$. Given an
element $\alpha\in H^2(\pi_1^{orb}(X), \mathbb S^1)$, we define
the twisted
orbifold K-theory of $X$ as $^\alpha K^*_{orb}(X) =
^\alpha K^*_{\pi}
(\tilde{M})$, where $\pi = \pi_1^{orb}(X)\times G$.
\end{definition}

We can also define the {\sl total twisted orbifold K-theory} of $X$ as

$$TK^*_{orb}(X) = \bigoplus_{\alpha\in H^2(\pi_1^{orb}(X),~ \mathbb S^1)}
~^\alpha K^*_{orb}(X).$$
This will also have a bigraded ring structure, and it will exhibit
an isomorphism between the $\alpha$ component and the $-\alpha$
component.

If $Y$, the orbifold universal cover of $X$, is actually a {\sl
manifold}, then $X$ is said to be a \emph{good} orbifold (see
\cite{MM}). The $G$ action on $\tilde{M}$ is free, and in this
case the $\alpha$--twisted orbifold cohomology will simply be
$^\alpha K^*_{\pi_1^{orb}(X)}(Y)$. 
For the case of a global quotient $X=Z/G$ and a class
$\alpha\in H^2(G,\mathbb S^1)$ it is not hard to verify
that in fact
$^{f*(\alpha)} K^*_{orb}(X)\cong ^\alpha K^*_G(Z)$, where
$f:\pi_1^{orb}(X)\to G$ is defined as before.

In the general case we note
that $\pi = \pi_1^{orb}(X)\times G$ acts on $\tilde{M}$ with
finite isotropy, hence we can make use of `twisted
Bredon cohomology' and a twisted version of the usual
Atiyah--Hirzebruch spectral sequence. Fix $\alpha\in
H^2(\pi_1^{orb}(X),\mathbb S^1)$, where $X$ is a compact orbifold.
There is a spectral sequence of the form
$$E_2=H^*_{\pi}(\tilde{M}; R_{\alpha}(-))\Rightarrow 
^\alpha K^*_{orb}(X)$$ The $E_1$ term will be a
chain complex built out of the twisted representation rings of the
stabilizers, all of which are finite. In many cases, this twisted
Atiyah--Hirzebruch spectral sequence will also collapse at $E_2$
after tensoring with the complex numbers. We believe that in fact
this must always be the case. In particular 
we conjecture that if 
$X$ denotes a compact good orbifold, with orbifold universal 
cover the manifold
$Y$, with $\Gamma = \pi_1^{orb}(X)$, and $\alpha\in H^2(\Gamma, \mathbb S^1)$,
then we have an additive decomposition
\[
^\alpha K^*_{\Gamma}(X)\otimes\mathbb C\cong \bigoplus_{(\gamma)}~
H^*(Hom_{Z_{\Gamma}(\gamma)}(C_*(Y^{\langle\gamma\rangle}), L_{\alpha}^{\gamma}))
\cong H^*_{orb}(X,\mathcal{L}_\alpha)
\]
where $(\gamma)$ ranges over conjugacy classes of elements of finite order
in $\Gamma$, $C_*(-)$ denotes the singular 
chains and $L_{\alpha}^{\gamma}$
is the character for $Z_{\Gamma}(\gamma)$ associated to the twisting.

\end{document}